\documentclass{article}
\usepackage{amssymb}
\usepackage{amsbsy}
\usepackage{amsmath}
\begin{document}
\begin{center}
{\LARGE \bf On Hilbert's sum type inequalities}
\end{center}

\vskip 15pt

\centerline{Chang-Jian Zhao\footnote{Research is supported by
National Natural Science Foundation of China (11371334).}}
\centerline{\it Department of Mathematics, China Jiliang
University, Hangzhou 310018, P.R.China} \centerline{\it Email:
chjzhao@163.com~~ chjzhao@cjlu.edu.cn}

\centerline{Wing-Sum Cheung\footnote{ Research is partially
supported by the Research Grants Council of the Hong Kong SAR,
China (Project No.: HKU7016/07P).}} \centerline{\it Department of
Mathematics, The University of Hong Kong, Pokfulam Road, Hong
Kong} \centerline{\it E-mail: wscheung@hku.hk}

\vskip 20pt

\begin{center}
\begin{minipage}{12cm}

{\large\bf Abstract}~ The main purpose of the present article is to give some new Hilbert's sum type inequalities, which in special cases yield the classical Hilbert's inequalities. Our results provide some new estimates to these types of inequalities.

{\large\bf MR (2010) Subject Classification} 26D15

\end{minipage}
\end{center}

\vskip 10pt \normalsize {\Large \bf 1 ~Introduction}\vskip 10pt

The classical Hilbert's inequalities can be stated as
follows ([1], p.226)

{\bf Hilbert's integral inequality}: {\it Let $f(x),g(x)\geq 0$, $0<\int_0^\infty
f^p(x)dx\leq \infty$ and $0<\int_0^\infty g^q(y)dy\leq \infty$. If
$p>1$ and $q=p/(p-1)$, then}
$$\int_0^\infty\int_0^\infty\frac{f(x)g(y)}{x+y}dxdy\leq\frac{\pi}{sin(\pi/p)}\left(\int_0^\infty
f^p(x)dx\right)^{1/p}\left(\int_0^\infty
g^q(y)dy\right)^{1/q}.\eqno(1.1)$$

{\bf Hilbert's discrete inequality}: {\it Let $a_{m}, b_{m}\geq 0$,  $0<\left(\sum_{1}^{\infty}a_{m}^{p}\right)^{1/p}\leq \infty$ and $0<\left(\sum_{1}^{\infty}a_{m}^{q}\right)^{1/q}\leq \infty$. If
$p>1$ and $q=p/(p-1)$, then}
$$\sum_{1}^{\infty}\sum_{1}^{\infty}\frac{a_{m}b_{m}}{m+n}\leq\frac{\pi}{\sin(\pi/p)}
\left(\sum_{1}^{\infty}a_{m}^{p}\right)^{1/p}\left(\sum_{1}^{\infty}b_{m}^{q}\right)^{1/q}.\eqno(1.2)$$

Hilbert's integral inequality and its discrete form were studied extensively and numerous
variants, generalizations, and extensions appeared in the
literature [2-12] and the references cited therein. The main purpose of the present article is to give some Hilbert's sum form inequalities, which in special cases yield
some Hilbert's type inequalities. Our main results are given in the following inequalities.

{\bf Hilbert's sum discrete inequality}~ {\it If $a_{m}\geq 0$, $b_{m}\geq 0$, $c_{m}\geq 0$, $d_{m}\geq 0$, where $m=0,1,\dots,$ and $p>1$ and $q=p/(p-1)$, and $0<\sum_{1}^{\infty}a_{m}^{p},\sum_{1}^{\infty}b_{m}^{q}, \sum_{0}^{\infty}c_{m}^{p},$ $\sum_{0}^{\infty}d_{m}^{q}<\infty$, then}
$$c_{0}d_{0}+\sum_{1}^{\infty}\sum_{1}^{\infty}\left(\frac{a_{m}b_{m}}{m+n}+\frac{c_{m}d_{m}}
{m+n+1}\right)\leq
\frac{\pi}{\sin(\pi/p)}\left(c_{0}^{p}+\sum_{1}^{\infty}\left(a_{m}^{p}+c_{m}^{p}\right)\right)^{1/p}
$$$$~~~~~~~~~~~~~~~~~~~~~~~~~~~~~~~~~~~~~~\times\left(d_{0}^{q}+\sum_{1}^{\infty}
\left(b_{m}^{q}+d_{m}^{q}\right)\right)^{1/q}.\eqno(1.3)$$

When $c_{m}=d_{m}=0$ and $m=0,1,\ldots,$ (1.3) becomes (1.2). When $a_{m}=b_{m}=0$ and $m=1,\ldots,$ (1.3) becomes (2.2) in the Section 2. Moreover, inequality (1.3) is just a special case of (2.5) established in the Section 2.

{\bf Hilbert's sum integral inequality}~ {\it If $p>1, \frac{1}{p}+\frac{1}{q}=1$, $f(x),g(y)>0, \lambda>2n, \gamma\in(n-\frac{\lambda}{p},\frac{\lambda}{q}-n),$ $n=0,1,\ldots,$ and $f(x)\geq 0$, the derivatives $f',f'',\ldots,f^{(n)}$ exists and positive and $f^{(n)}\in L(0,\infty), (n=0,12,\ldots),(f^{(0)}:=f)$, where $L^{p}$ denotes the space of all Lebesgue integrable function (If $p=1$, we obtain $L$), and let $f(0)=f'(0)=f''(0)=\cdots=f^{(n-1)}(0)=0$ such that $\int_{0}^{\infty}x^{p(n+1)-p\gamma-\lambda-1}(f^{(n)}(x))^{p}dx<\infty$. Moreover, $g(y)$ and $f(x)$ have the same conditions as above, then
$$\int_0^\infty\int_0^\infty\frac{f(x)g(y)}{(x+y)^{\lambda}}dxdy\leq\frac{1}{2}
\left(\int_{0}^{\infty}x^{p(n+1)-\lambda-p\gamma-1}(Cx^{p\gamma}+C')(f^{(n)}(x))^p dx\right)^{1/p}$$
$$~~~~~~~~~~~~~~~~~~~~~~~~~~\times\left(\int_{0}^{\infty}y^{q(n+1)-\lambda-q\gamma-1}(Cy^{q\gamma}+C')(g^{(n)}(y))^p dy\right)^{1/q},\eqno(1.4)$$
where}
$$C=\frac{\Gamma(\frac{\lambda}{p}-n)\Gamma(\frac{\lambda}{q}-n)}{\Gamma(\lambda)}~~~
{\it and}~~~
C'=\frac{\Gamma(\frac{\lambda}{p}-\gamma-n)\Gamma(\frac{\lambda}{q}-\gamma-n)}
{\Gamma(\lambda)}.\eqno(1.5)$$

Inequality (1.4) is just a special case of (2.7) established in the Section 2.

\vskip 10pt \normalsize {\Large \bf 2 ~Main results}\vskip 10pt

{\bf Lemma 2.1}~ ([13, p.26] p.39) {\it If $a_{i}\geq 0$,
$b_{i}>0, i=1,\ldots,m$ and $\sum_{i=1}^{m}\alpha_{i}=1$,
then
$$\left(\prod_{i=1}^{m}(a_{i}+b_{i})\right)^{\alpha_{i}}\geq\left(\prod_{i=1}^{m}a_{i}\right)^{\alpha_{i}}
+\left(\prod_{i=1}^{m}b_{i}\right)^{\alpha_{i}},\eqno(2.1)$$
with equality if and only if} $a_{1}/b_{1}=\cdots=a_{m}/b_{m}.$

{\bf Lemma 2.2}~ ([1]) {\it If $c_{m}\geq 0$, $d_{m}\geq 0$, where $m=0,1,\dots,$ and $p>1$ and $q=p/(p-1)$, and $0<\sum_{0}^{\infty}c_{m}^{p},$ $\sum_{0}^{\infty}d_{m}^{q}<\infty$, then}
$$\sum_{0}^{\infty}\sum_{0}^{\infty}\frac{c_{m}d_{m}}{m+n+1}\leq\frac{\pi}{\sin(\pi/p)}
\left(\sum_{0}^{\infty}c_{m}^{p}\right)^{1/p}\left(\sum_{0}^{\infty}d_{m}^{q}\right)^{1/q}.\eqno(2.2)$$

{\bf Lemma 2.3}~ ([3]) {\it If $p>1$ and $q=p/(p-1)$, and $f(x)>0$, the derivatives $f',f'',\ldots,f^{(n)}$ exists and positive and $f^{(n)}\in L(0,\infty), (n=0,12,\ldots),(f^{(0)}:=f)$, where $L^{p}$ denotes the space of all Lebesgue integrable function (If $p=1$, we obtain $L$), and let $f(0)=f'(0)=f''(0)=\cdots=f^{(n-1)}(0)=0$ such that $\int_{0}^{\infty}x^{p(n+1)-\lambda-1}(f^{(n)}(x))^{p}dx<\infty$. Moreover, $g(y)$ and $f(x)$ have the same conditions as above, then
$$\int_0^\infty\int_0^\infty\frac{f(x)g(y)}{(x+y)^{\lambda}}dxdy\leq C\left(\int_0^\infty
x^{p(n+1)-\lambda-1}(f^{(n)}(x))^p dx\right)^{1/p}$$
$$~~~~~~~~~~~~~~~~~~~~~~~~~\times\left(\int_0^\infty y^{q(n+1)-\lambda-1}(g^{(n)}(y))^q dy\right)^{1/q}，\eqno(2.3)$$
where $C$ is as in (1.5), and which is the best possible. Here $\Gamma(u)$ is the Gamma function}.

{\bf Lemma 2.4}~ ([12]) {\it Let $p>1$ and $q=p/(p-1)$, $\lambda>2n, \gamma\in(n-\frac{\lambda}{p},\frac{\lambda}{q}-n), n=0,1,\ldots$ and $f(x)>0$, the derivatives $f',f'',\ldots,f^{(n)}$ exists and positive and $f^{(n)}\in L(0,\infty), (n=0,12,\ldots),(f^{(0)}\\:=f)$, where $L^{p}$ denotes the space of all Lebesgue integrable function (If $p=1$, we obtain $L$), and let $f(0)=f'(0)=f''(0)=\cdots=f^{(n-1)}(0)=0$ such that $\int_{0}^{\infty}x^{p(n+1)-p\gamma-\lambda-1}(f^{(n)}(x))^{p}dx<\infty$. Moreover, $g(y)$ and $f(x)$ have the same conditions as above, then
$$\int_0^\infty\int_0^\infty\frac{f(x)g(y)}{(x+y)^{\lambda}}dxdy\leq C'\left(\int_0^\infty
x^{p(n+1)-p\gamma-\lambda-1}(f^{(n)}(x))^p dx\right)^{1/p}$$
$$~~~~~~~~~~~~~~~~~~~~~~~~~~\times\left(\int_0^\infty y^{q(n+1)-p\gamma-\lambda-1}(g^{(n)}(y))^q dy\right)^{1/q}.\eqno(2.4)$$
where $C$ is as in (1.5), and which is the best possible}.

{\bf Theorem 2.1}~ {\it If $a_{m}\geq 0$, $b_{m}\geq 0$, $c_{m}\geq 0$, $d_{m}\geq 0$, where $m=0,1,\dots,$ and $k\in{\Bbb N}^{+}$ and $p>1$ and $q=p/(p-1)$, and $\sum_{1}^{\infty}a_{m}^{p},\sum_{1}^{\infty}b_{m}^{q}, \sum_{0}^{\infty}c_{m}^{p},$ $\sum_{0}^{\infty}d_{m}^{q}<\infty$, then}
$$kc_{0}d_{0}+\sum_{1}^{\infty}\sum_{1}^{\infty}\left(\frac{a_{m}b_{m}}{m+n}+\frac{kc_{m}d_{m}}
{m+n+1}\right)\leq
\frac{\pi}{\sin(\pi/p)}\left(kc_{0}^{p}+\sum_{1}^{\infty}\left(a_{m}^{p}+kc_{m}^{p}\right)\right)^{1/p}
$$$$~~~~~~~~~~~~~~~~~~~~~~~~~~~~~~~~~~~~~~~~\times\left(kd_{0}^{q}+\sum_{1}^{\infty}
\left(b_{m}^{q}+kd_{m}^{q}\right)\right)^{1/q}.\eqno(2.5)$$

{\it Proof}~ First, we prove that (2.5) holds for $k=1$. Noting that
$$\sum_{1}^{\infty}\sum_{1}^{\infty}\frac{a_{m}b_{m}}{m+n}\leq\frac{\pi}{\sin(\pi/p)}
\left(\sum_{1}^{\infty}a_{m}^{p}\right)^{1/p}\left(\sum_{1}^{\infty}b_{m}^{q}\right)^{1/q},$$
and
$$\sum_{0}^{\infty}\sum_{0}^{\infty}\frac{c_{m}d_{m}}{m+n+1}\leq\frac{\pi}{\sin(\pi/p)}
\left(\sum_{0}^{\infty}c_{m}^{p}\right)^{1/p}\left(\sum_{0}^{\infty}d_{m}^{q}\right)^{1/q}.$$
Hence
\begin{eqnarray*}c_{0}d_{0}+\sum_{1}^{\infty}\sum_{1}^{\infty}\left(\frac{a_{m}b_{m}}{m+n}+\frac{c_{m}d_{m}}
{m+n+1}\right)&\leq&\frac{\pi}{\sin(\pi/p)}\Bigg(\left(\sum_{1}^{\infty}a_{m}^{p}\right)^{1/p}\left(\sum_{1}^{\infty}b_{m}^{q}\right)^{1/q}\\
&+&\left(\sum_{0}^{\infty}c_{m}^{p}\right)^{1/p}\left(\sum_{0}^{\infty}d_{m}^{q}\right)^{1/q}\Bigg)\\
&\leq&\frac{\pi}{\sin(\pi/p)}\left(c_{0}^{p}+\sum_{1}^{\infty}\left(a_{m}^{p}+c_{m}^{p}\right)\right)^{1/p}
\left(d_{0}^{q}+\sum_{1}^{\infty}\left(b_{m}^{q}+d_{m}^{q}\right)\right)^{1/q}.\end{eqnarray*}
This shows (2.5) right for $k=1$.

Suppose that (2.5) holds when $k=r-1$, we have
$$(r-1)c_{0}d_{0}+\sum_{1}^{\infty}\sum_{1}^{\infty}\frac{a_{m}b_{m}+(r-1)c_{m}d_{m}}{m+n+1}\leq
\frac{\pi}{\sin(\pi/p)}\left((r-1)c_{0}^{p}+\sum_{1}^{\infty}\left(a_{m}^{p}+(r-1)c_{m}^{p}\right)\right)^{1/p}
$$$$~~~~~~~~~~~~~~~~~~~~~~~~~~~~~~~~~~~~~~~~~~~~~\times\left((r-1)d_{0}^{q}+\sum_{1}^{\infty}\left(b_{m}^{q}+(r-1)
d_{m}^{q}\right)\right)^{1/q}.\eqno(2.6)$$
From (2.2) and (2.6), we have
\begin{eqnarray*}rc_{0}d_{0}+\sum_{1}^{\infty}\sum_{1}^{\infty}\frac{a_{m}b_{m}+rc_{m}d_{m}}{m+n+1}
&\leq&\frac{\pi}{\sin(\pi/p)}\Bigg\{\left((r-1)c_{0}^{p}+\sum_{1}^{\infty}\left(a_{m}^{p}+(r-1)c_{m}^{p}\
\right)\right)^{1/p}\\
&\times&\left((r-1)d_{0}^{q}+\sum_{1}^{\infty}\left(b_{m}^{q}+(r-1)
d_{m}^{q}\right)\right)^{1/q}\\
&+&\left(\sum_{0}^{\infty}c_{m}^{p}\right)^{1/p}\left(\sum_{0}^{\infty}d_{m}^{q}\right)^{1/q}\Bigg\}\\
&\leq&\frac{\pi}{\sin(\pi/p)}\left(rc_{0}^{p}+\sum_{1}^{\infty}\left(a_{m}^{p}+rc_{m}^{p}\right)\right)^{1/p}\\
&\times&\left(rd_{0}^{q}+\sum_{1}^{\infty}\left(b_{m}^{q}+rd_{m}^{q}\right)\right)^{1/q}.\end{eqnarray*}
This shows that (2.5) is correct if $m=r-1$, then $m=r$ is also correct. Hence (2.5) is right for any $m\in{\Bbb N}^{+}$.

This proof is complete.\hfill $\Box$

{\bf Theorem 2.2}~ {\it Let $p,q$, $\gamma$, $\lambda$, $f(x),g(y)$ are as in} (2.4). {\it If $m\in{\Bbb N}^{+}$, then
$$\int_0^\infty\int_0^\infty\frac{f(x)g(y)}{(x+y)^{\lambda}}dxdy\leq\frac{1}{m+1}\left(\int_{0}^{\infty}x^{p(n+1)-\lambda-p\gamma-1}(Cx^{p\gamma}+mC')(f^{(n)}(x))^p dx\right)^{1/p}$$
$$~~~~~~~~~~~~~~~~~~~~\times\left(\int_{0}^{\infty}y^{q(n+1)-\lambda-q\gamma-1}(Cy^{q\gamma}+mC')(g^{(n)}(y))^p dy\right)^{1/q},\eqno(2.7)$$
where $C'$ and $C$ are as in} (1.5).

{\it Proof}~ First, we prove that (2.7) holds for $m=1$. From Lemmas 2.3 and 2.4, we hgave
$$\int_0^\infty\int_0^\infty\frac{f(x)g(y)}{(x+y)^{\lambda}}dxdy\leq C\left(\int_0^\infty
x^{p(n+1)-\lambda-1}(f^{(n)}(x))^p dx\right)^{1/p}$$
$$~~~~~~~~~~~~~~~~~~~~~~~~~~~~\times\left(\int_0^\infty y^{q(n+1)-\lambda-1}(g^{(n)}(y))^q dy\right)^{1/q}.\eqno(2.8)$$
where $C$ is as in (1.5)
$$\int_0^\infty\int_0^\infty\frac{f(x)g(y)}{(x+y)^{\lambda}}dxdy\leq C'\left(\int_0^\infty
x^{p(n+1)-p\gamma-\lambda-1}(f^{(n)}(x))^p dx\right)^{1/p}$$
$$~~~~~~~~~~~~~~~~~~~~~~~~~~~\times\left(\int_0^\infty y^{q(n+1)-p\gamma-\lambda-1}(g^{(n)}(y))^q dy\right)^{1/q}.\eqno(2.9)$$
where $C'$ is as in (1.5)

From (2.8) and (2.9), and in view of Lemma 2.1, we obtain
\begin{eqnarray*}\int_0^\infty\int_0^\infty\frac{f(x)g(y)}{(x+y)^{\lambda}}dxdy&\leq& \frac{C}{2}\left(\int_0^\infty x^{p(n+1)-\lambda-1}(f^{(n)}(x))^p dx\right)^{1/p}\left(\int_0^\infty
y^{q(n+1)-\lambda-1}(g^{(n)}(y))^q dy\right)^{1/q}\\
&+&\frac{C'}{2}\left(\int_0^\infty
x^{p(n+1)-p\gamma-\lambda-1}(f^{(n)}(x))^p dx\right)^{1/p}\\
&\times&\left(\int_0^\infty y^{q(n+1)-q\gamma-\lambda-1}(g^{(n)}(y))^q dy\right)^{1/q}\\
&\leq&\frac{1}{2}\left(\int_{0}^{\infty}x^{p(n+1)-\lambda-p\gamma-1}(Cx^{p\gamma}+C')(f^{(n)}(x))^p dx\right)^{1/p}\\
&\times&\left(\int_{0}^{\infty}y^{q(n+1)-\lambda-q\gamma-1}(Cy^{q\gamma}+C')(g^{(n)}(y))^p dy\right)^{1/q}.\end{eqnarray*}

This shows (2.7) right for $m=1$.

Suppose that (2.7) holds when $m=r-1$, we have
$$\int_0^\infty\int_0^\infty\frac{f(x)g(y)}{(x+y)^{\lambda}}dxdy\leq\frac{1}{r}\left(\int_{0}^{\infty}x^{p(n+1)-\lambda-p\gamma-1}(Cx^{p\gamma}+(r-1)C')(f^{(n)}(x))^p dx\right)^{1/p}$$
$$~~~~~~~~~~~~~~~~~~~~~~~~~~~~~~~~~~~\times\left(\int_{0}^{\infty}y^{q(n+1)-\lambda-q\gamma-1}(Cy^{q\gamma}+(r-1)C')(g^{(n)}(y))^p dy\right)^{1/q}.\eqno(2.10)$$

From (2.4), (2.10) and by using
\begin{eqnarray*}\int_0^\infty\int_0^\infty\frac{(r+1)f(x)g(y)}{(x+y)^{\lambda}}dxdy&\leq&
\Bigg(\left(\int_{0}^{\infty}x^{p(n+1)-\lambda-p\gamma-1}(Cx^{p\gamma}+(r-1)C')(f^{(n)}(x))^p dx\right)^{1/p}\\
&\times&\left(\int_{0}^{\infty}y^{q(n+1)-\lambda-q\gamma-1}(Cy^{q\gamma}+(r-1)C')(g^{(n)}(y))^p dy\right)^{1/q}\\
&+&C'\left(\int_0^\infty x^{p(n+1)-p\gamma-\lambda-1}(f^{(n)}(x))^p dx\right)^{1/p}\\
&\times&\left(\int_0^\infty y^{q(n+1)-p\gamma-\lambda-1}(g^{(n)}(y))^q dy\right)^{1/q}\Bigg)\\
&\leq&\left(\int_{0}^{\infty}x^{p(n+1)-\lambda-p\gamma-1}(Cx^{p\gamma}+rC')(f^{(n)}(x))^p dx\right)^{1/p}\\
&\times&\left(\int_{0}^{\infty}y^{q(n+1)-\lambda-q\gamma-1}(Cy^{q\gamma}+rC')(g^{(n)}(y))^p dy\right)^{1/q}.\end{eqnarray*}

This shows that (2.7) is correct if $m=r-1$, then $m=r$ is also correct. Hence (2.7) is right for any $n\in{\Bbb N}$.

This completes the proof of Theorem 2.2.\hfill $\Box$

\vskip 8pt {\small

\end{document}